\documentclass[makeidx]{amsart}
\usepackage[dvips]{graphicx}
\usepackage{a4wide}
\usepackage{latexsym}
\usepackage{amsmath}
\usepackage{amssymb}
\usepackage{amsthm}
\usepackage{amsbsy}

\newtheorem{theorem}{Theorem}[section]
\newtheorem{lemma}[theorem]{Lemma}
\newtheorem{e-proposition}[theorem]{Proposition}

\newtheorem{e-definition}[theorem]{Definition\rm}
\newtheorem{remark}{\it Remark\/}

\newcommand{\R}{\mathbb{R}}
\newcommand{\C}{\mathbb{C}}
\newcommand{\F}{\mathcal{F}}

\begin{document}

\title[The Reeb foliation arises as a family of Legendrian submanifolds]
{
The Reeb foliation\\ arises as a family of Legendrian submanifolds\\
at the end of a deformation of the standard $S^3$ in $S^5$
}
\author[A.~Mori]
{
Atsuhide MORI
}

\address{Graduate~School~of~Science, Osaka~University, 
Toyonaka, Osaka~560-0043, Japan}

\email{ka-mori@ares.eonet.ne.jp}



\begin{abstract}
We realize the Reeb foliation of $S^3$ as a family 
of Legendrian submanifolds of the unit $S^5\subset \C^3$,
Moreover we construct a deformation of the 
standard contact $S^3$ in $S^5$,
via a family of contact submanifolds,
into this realization.
\end{abstract}

\maketitle

\section{Introduction}
The Reeb foliation is a codimension one 
smooth foliation of the $3$-sphere $S^3$ 
obtained by glueing two Reeb components
$S^1\times D^2$ and $D^2\times S^1$.
Since the one-sided holonomies of the Reeb components
along $\{1\}\times \partial D^2$
and $\partial D^2\times \{1\}$
are trivial, the Reeb foliation
is not analytic (``Haefliger's remark'').

On the other hand the $1$-jet space $J^1(\R^n,\R)\approx \R^{2n+1}$ 
for a function of $n$ variables carries the canonical contact structure.
It is contactomorphic to the unit sphere $S^{2n+1}\subset \C^{n+1}$ minus 
any point. Here $S^{2n+1}$ has the standard contact form
$\alpha=\sum_{i=1}^{n+1} r_i^2d\theta_i|S^{2n+1}$
($r_i=|z_i|$, $\theta_i=\arg z_i$ for coordinates $z_i$ of $\C^{n+1}$).
Thus we may regard a codimension-$n$ submanifold $M^{n+1}\subset S^{2n+1}$ 
as a system of $n$ first-order partial differential equations 
(for implicit functions).
If $\alpha\wedge d\alpha|M^{n+1}=0$ and $\alpha|M^{n+1}\ne 0$,
the system is completely integrable and regular, and therefore
defines a codimension one foliation $\F$ on $M^{n+1}$.
The leaves of $\F$ are Legendrian
submanifolds of $S^{2n+1}$ corresponding to the solutions.

In this article we construct an embedding of $S^3$ into the standard $S^5$
so that the image has the Reeb foliation $\F$ by
Legendrian submanifolds. This example shows that even 
a non-taut foliation can be a family of Legendrian submanifolds
of $J^1(\R^n,\R)$. Moreover we prove
\begin{theorem}
\label{Main}
There exists a smooth family $\{M^3_t\}_{t\in[0,3/2)}$ of 
codimension-$2$ submanifolds of $S^5$ 
such that 
\begin{enumerate}
\item[(1)] $M_0^3$ is the standard $S^3(\subset\C^2\subset\C^3)$,
\item[(2)] $M_t^3$ is an embedded contact submanifold for $0\le t<1$,
\item[(3)] $M_1^3$ admits a Reeb foliation by 
injectively immersed Legendrian submanifolds of $S^5$, and
\item[(4)] $-M_t^3$ is an embedded 
overtwisted contact submanifold for $1<t< 3/2$.
\end{enumerate}
\end{theorem}
The foliated submanifold $M^3_1$ is obtained by
joining two great circles $\{r_1=1\},\, \{r_2=1\}\subset S^5$ 
through the Legendrian torus
$T=\{r_1=r_2=r_3=1/\sqrt{3},\, \theta_1+\theta_2+\theta_3=0\}$.
The family $M_t^3$ is obtained as a byproduct in the process
of isotoping $M_1\subset S^5$ to the unknot. 
The author is seeking the converse approach,
i.e., to find a foliated submanifold by using 
contact topology or open-books
(see Remark 1 in \S 2).
\section{Proof and remark}
\begin{proof}
Let $\pi$ be the natural projection of $S^5$ to 
the $2$-simplex $\Delta=\{(r_1^2,r_2^2,r_3^2)\, |\, r_1^2+r_2^2+r_3^2=1\}\subset\R^3$,
which sends the Legendrian $2$-torus 
$T=\{r_1=r_2=r_3=1/\sqrt{3},\, \theta_1+\theta_2+\theta_3=0\}\subset S^5$
to the barycenter $G$.
The set $\Gamma=\pi^{-1}(\partial\Delta)$ contains the great circles 
$\pi^{-1}(\{V_1,V_2,V_3\})$ where $V_i$ denotes the vertex $r_i^2=1$.
Except them $\pi|\Gamma$ is a $T^2$-fibration.
On the other hand, $\pi|(S^5\setminus \Gamma)$ 
is a $T^3$-fibration. 
Now we take the coordinates $(x,y)$ on $\Delta$ 
by putting $\displaystyle\overrightarrow{OP\,\,}
=\overrightarrow{OG\,\,}+x\overrightarrow{GV_1\,}
+y\overrightarrow{GV_2\,}$ for $P\in\Delta$, i.e.,
$$
3r_1^2=1+2x-y(\ge0),\quad 3r_2^2=1-x+2y(\ge 0),\quad
\textrm{and}\quad
3r_3^2=1-x-y(\ge 0).
$$
Let $M^3_0$ be the standard $S^3=\pi^{-1}(\overline{V_1 V_2})$.
We deform $M^3_0$ with the help of a certain family of simple curves 
$C_t: x=x_t(s),\, y=y_t(s),\, 
-\delta\le s\le \delta$ 
depicted in Fig.\ref{fig} ($0<\delta\ll1,\,0\le t\le 3/2$). 
Note that $C_1$ has a break point $G$ 
while $x_1(s)$ and $y_1(s)$ are smooth on $(-\delta,\delta)$.

\begin{figure}[h]
\centering
\includegraphics[height=35mm]{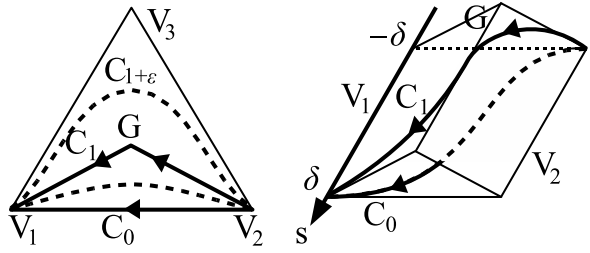}
\caption{The curve $C_t$ on $\Delta$ and its parametrization by $s$}
\label{fig}
\end{figure}

\noindent
We generate $M_t^3\subset S^5$ 
by moving the intersection of the ``wall''
$W_s=\textrm{cl}\{\theta_1+\theta_2+\theta_3=s\}\subset S^5$ with
the fibre $\pi^{-1}(x_t(s),y_t(s))$ 
for $-\delta\le s\le \delta$.
Then we can see that $M_t^3$ realizes the join of
two large circles $\pi^{-1}(V_2)$ and $\pi^{-1}(V_1)$.
Now we give a precise definition of the curve $C_t$.
Put $\displaystyle \varphi_0(u)=\frac{1}{2}(1+u)$ for $u\in [-1,1]$, and
take a smooth function $\varphi_1(u)$ and a 
smooth odd function $s(u)$ such that
$$
\begin{array}{l}
\varphi_1(u)=0\quad(-1\le u\le 0),
\quad
\varphi_1'(u)>0\quad(0< u\le 1),
\quad
\varphi_1(u)=\varphi_0(u)\quad(1/2\le u\le 1),
\\
s'(u)>0\quad(-1<u<1),
\quad
s(1)=\delta,\quad s(-1)=-\delta,
\quad\textrm{and}\quad
s(u) \textrm{ is } C^\infty \textrm{-tangent to } \pm \delta.
\nonumber
\end{array}
$$
The inverse function 
$u(s)$ of $s(u)$ is defined on $[-\delta,\delta]$.
It is smooth on $(-\delta,\delta)$ ($u'(\pm\delta)=+\infty$).
We put
$\varphi_t(u)=(1-t)\varphi_0(u)+t\varphi_1(u)$,
and take the curve
$$
C_t:\quad x=x_t(s)=\varphi_t(u(s)),\quad y=y_t(s)=\varphi_t(u(-s)),
\quad -\delta \le s \le \delta.
$$

Next we show that $M_t^3$ is a smooth submanifold.
By moving the $2$-torus $(M^3_t\setminus \Gamma)\cap W_s$
for $-\delta<s<\delta$,
we see that $M^3_t\setminus \Gamma$ is diffeomorphic to 
$T^2\times(-\delta,\delta)$. 
Moreover $M_t^3$ is topologically the join $S^1\star S^1\approx S^3$. 
Thus it only remains for us 
to examine the smoothness of $M_t^3$ along $M_t^3\cap\Gamma$. 
We restrict ourself 
to the connected component of $M_t^3\cap\Gamma$ 
corresponding to $s=+\delta$ and omit the other component. 
We put
$$
\widetilde{M_t^3}:\left\{
\begin{array}{l}
\displaystyle r_1^2+r_2^2+r_3^2=1\\
\displaystyle 3r_2^2=1 -\frac{1}{2}(1+u) +(1-t)(1-u)\\
\displaystyle 3r_3^2=1 -\frac{1}{2}(1+u) -\frac{1-t}{2}(1-u)
=\frac{t}{3-2t}\cdot 3r_2^2\\
\displaystyle \theta_1+\theta_2+\theta_3=1
\end{array}\right.
$$
where $u\in [1/2, 1]$ is a parameter to be eliminated.
Then $\{\theta_1=\mathrm{const}\}\subset \widetilde{M_t^3}$ 
is a smooth disk 
since it tangents to the real $2$-plane
$\displaystyle \left\{z_1=\exp{\sqrt{-1}\theta_1},\quad z_3=
\overline{z_2}\cdot \sqrt{\frac{t}{3-2t}}\exp{\{\sqrt{-1}(1-\theta_1)\}}
\right\}\subset \C^3$ at $u=1$.
Since the function $s(u)$ smoothly tangents to $\delta$ at $u=1$, $M_t^3$ is a smooth $3$-sphere.

Next we consider the (non-)integrability of the restriction $\lambda_t=\alpha|
M_t^3$ of the standard contact form 
$\alpha=r_1^2d\theta_1+r_2^2d\theta_2+r_3^2d\theta_3|S^5$.
Using $(\theta_1,\theta_2,s)$ as coordinates of $M_t^3\setminus\Gamma$, 
we can write
$$\lambda_t=x_t(s)d\theta_1+y_t(s)d\theta_2+(1-x_t(s)-y_t(s))ds.$$
Here the sign of $\lambda_t\wedge d\lambda_t$ 
with respect to $d\theta_1 \wedge d\theta_2 \wedge ds>0$
coincides with that of $x_t'(s)y_t(s)-x_t(s)y_t'(s)$, 
and that of $1-t$. More generally, 
if a submanifold $M^3(\approx T^2\times \R)\subset S^5$
is presented by a simple curve 
$C:\, x=x(s), \, y=y(s)$ on $\textrm{int}\Delta$, 
the nagative areal verocity $x'(s)y(s)-x(s)y'(s)$ 
still presents the non-integrability of $\alpha|M^3$. 
In the case where $t=1$, the integrability means the vanishing
of the areal verocity. That is why the curve $C_1$ is
broken into two rays to/from the origin $G$, 
and $M_1^3$ is non-analytic. 

On the other hand, for cylindrical 
coordinates $(\theta_1,(r_2,\theta_2))$, 
$\mu_t=\alpha|\widetilde{M_t^3}$ and $\mu_t\wedge d\mu_t$ are
written as
$$
\mu_t=\left(1-\frac{3}{3-2t}r_2^2\right)d\theta_1
+\frac{3(1-t)}{3-2t}r_2^2 d\theta_2
\quad\textrm{and}\quad
\mu_t\wedge d\mu_t=\frac{6(1-t)}{3-2t}d\theta_1\wedge(r_2dr_2\wedge d\theta_2).
$$ 
This implies that the sign of $\lambda_t\wedge d\lambda_t$ everywhere coincides with that of $1-t$.

Now we show that the foliation of $M^3_1$ is a Reeb foliation. 
The definition of $M^3_1$ is
$$
\left\{
\begin{array}{l}
\displaystyle 3r_1^2=1+2\varphi(u(s))-\varphi(u(-s))
\\
\displaystyle 3r_2^2=1-\varphi(u(s))+2\varphi(u(-s))
\\
\displaystyle 3r_3^2=1-\varphi(u(s))-\varphi(u(-s))
\\
\theta_1+\theta_2+\theta_3=s
\end{array}
\right.
$$
where $s\in[-\delta,\delta]$ is a parameter to be eliminated. 
On the open solid torus $H=\{s>0\}\subset M^3_1$, we have
$$
\alpha|H=\varphi(u(s))d\theta_1+\{1-\varphi(u(s))\}ds.
$$
Thus the surface of $\theta_2$-revolution of the graph of 
$\displaystyle\theta_1=\int \frac{\varphi(u(s))-1}{\varphi(u(s))}ds$
is a leaf. Similarly, we can describe the foliation on $\{s<0\}$.
These foliations spiral into $T$ and form a 
transversely oriented Reeb foliation,
to which the positive Hopf link $\{r_1=1\}\cup\{r_2=1\}$ 
is positively transverse.

Finally we see from $d(\theta_1+\theta_2)\wedge d\lambda_t
=\{x_t'(s)-y_t'(s)\}d\theta_1\wedge d\theta_2\wedge ds>0\,\, (t\ne 1)$
that the positive Hopf band $\ker(d\theta_1+d\theta_2)$ 
is a supporting open-book for $0\le t<1$. 
On the other hand, the negative Hopf band
$\ker(-d\theta_1-d\theta_2)$ on $-M_t(\approx S^3)$
is a supporting open-book for $1< t<3/2$. 
Thus $-M_t^3$ is overtwisted. 
Indeed it has the half-Lutz tube $\{x_t(s)\le 0\}$. 
Moreover, since we can reverse the orientation of $S^3$ by a diffeotopy, 
we obtain the following ``negative stabilization'' lemma.
This ends the proof.
\end{proof}
\begin{lemma}
The overtwisted contact submanifold 
$-M_{5/4}^3\subset S^5$ is 
diffeotopic to the standard $S^3\subset S^5$. Particularly
$-M_{5/4}^3$ is differential topologically unknotted, but contact topologically knotted.
\end{lemma}

\begin{remark}
Any closed oriented $3$-manifold admits an open-book decomposition
(Alexander \cite{Alexander}). We can associate to it a contact structure
(Thurston-Winkelnkemper \cite{TW}) as well as a spinnable foliation (see \cite{Mori1}). 
Further any contact structure is supported by an open-book decomposition
(Giroux \cite{Giroux}). Using this result, the author constructed 
a certain immersion of any contact $3$-manifold 
into $J^1(\R^2,\R)$ or $S^5$ (\cite{Mori2}). This construction was
generalized to any dimension, i.e., $M^{2n+1}\to J^1(\R^{2n},\R)$ or $S^{4n+1}$
by Mart\'inez Torres (\cite{Martinez}).
The author proved that any/some contact structure of $M^3$
can be deformed into some/any spinnable foliation
(\cite{Mori1}, see also \cite{Etnyre}). He also proved that a certain
higher dimensional contact structure can be deformed into a foliation (\cite{Mori3}).
It is interesting to generalize the present result to these cases.
\end{remark}

\section*{Acknowledgements}
The author would like to thank the anonymous reviewer(s) for 
encouragement to include intuitive descriptions and a figure,
which have made the article easier to access for general readers.

\end{document}